# TOPOLOGICAL ENTROPY OF GENERALIZED POLYGON EXCHANGES

EUGENE GUTKIN AND NICOLAI HAYDN

ABSTRACT. We obtain geometric upper bounds on the topological entropy of generalized polygon exchange transformations. As an application of our results, we show that billiards in polygons and rational polytops have zero topological entropy.

## 1. INTRODUCTION

Let $T : M \mapsto M$ be a homeomorphism of a compact manifold. The topological entropy is an important invariant that characterizes by a number, $h(T) \geq 0$, the complexity of the orbits of $T$. There is a lot of information on $h(T)$ in the literature. The results depend on the geometric context for $T : M \mapsto M$. Let us mention two types of theorems on topological entropy.

Suppose that $T$ comes from the geodesic flow on a compact Riemannian manifold of negative curvature. Then one has explicit bounds on $h(T)$ from above and from below in terms of the curvature of the manifold, its volume, diameter, etc. (see, e.g., [18, 9, 14, 15]).

Theorems of the other type provide estimates on $h(T)$ in terms of the growth rates of various geometric quantities on $M$, under the iterations of $T$ (see, e.g., [17, 20, 16, 5]). A recent result that is particularly relevant for our direction belongs to S. Newhouse ([16], see [17] for a closely related theorem). Theorem 1 below is a special case of Newhouse's result.

Assume that $M$ is a compact $C^1$ manifold and that $T : M \mapsto M$ is a $C^1$ diffeomorphism. Let $\gamma$ be a compact $C^1$ curve on $M$, and denote by $|\gamma|$ the length of $\gamma$ with respect to a Riemannian metric on $M$. The number

$$g_1(T) = \sup_\gamma \limsup_{n \to \infty} \frac{1}{n} \log |T^n \gamma|$$

does not depend on the choice of Riemannian metric, and we say that $g_1(T)$ is the *growth rate of the length of curves* on $M$.

**Theorem 1.** *Let $M$ be a compact $C^2$ surface, and let $T : M \mapsto M$ be a $C^{1+\epsilon}$ diffeomorphism. Then $h(T) \leq g_1(T)$.*

A metric entropy version of this theorem goes back to [12].

Theorem 1 and other results on $h(T)$ require that $T : M \mapsto M$ be sufficiently smooth. There is no literature on the topological entropy of transformations with

Received by the editors June 8, 1994.
1991 *Mathematics Subject Classification*. Primary 28D20, 58F11, 52B99.
Both authors were partially supported by the NSF.

singularities, despite the fact that they naturally arise in dynamics, e.g., in billiard flows [10, 11, 6, 8]. In the present work we introduce a class of transformations with controlled singularities: the *generalized polyhedron exchanges*; see Definitions 2 and 3. These mappings are the analogs, in many dimensions, of the interval exchange transformations [3]. The dimension of the polyhedra in question is arbitrary, but our techniques work especially well in two dimensions, for the *generalized polygon exchanges* (GPEs). The framework of GPEs is wide enough to include transformations arising in very different contexts. For instance, the Baker transformation [1], on one hand, and the polygonal billiard and outer billiard maps [6, 8], on the other hand, are GPEs.

The *topological entropy*, $h(T) \geq 0$, of a generalized polyhedron exchange characterizes the exponential complexity of regular orbits (§3). Our main result (Theorem 5) is an analog of Theorem 1 for the GPEs. It says that $h(T)$ is bounded above by the growth rate of the length of the singular set of $T^n$ as $n \to \infty$.

The *Lyapunov exponent*, $\lambda(x)$, is defined on the set of regular points for a generalized polyhedron exchange. It also provides an upper bound on the topological entropy of a GPE, $h(T) \leq \sup_x \lambda(x)$ (under some technical assumptions, see Theorem 6).

To illustrate the main result, we present three applications of Theorem 5 in §4. Our first application is to the Euclidean polygon exchanges, which are a direct analog, in two dimensions, of the interval exchange transformations (see the examples in §2). By Theorem 7, the topological entropy of any Euclidean polygon exchange is zero.

Then we give two applications of Theorem 5 to billiards. Theorem 8 states that the topological entropy of any polygonal billiard is zero. This implies the subexponential growth for the cardinality of the defining partition and for other geometric quantities associated with a polygonal billiard [10]. The earlier proofs of the subexponential growth (for polygonal billiards) involved a detailed analysis of the dynamics and geometry of polygons [10, 4]. Our proof is based on the observation that the length of the singular set for $T^n$ grows at most quadratically, hence Theorem 5 applies. Our third application is to the *directional billiards in a rational polytop* (Theorem 9). The corresponding Poincaré map, $T$, is a GPE. The length of the singular set grows subexponentially, hence, by Theorem 5, $h(T) = 0$.

## 2. Generalized polygon exchanges

By a *polyhedron* $P$, $\dim P = n$, we mean a compact Euclidean polyhedron in $\mathbf{R}^n$, $\operatorname{int} P \neq \emptyset$. The $m$-dimensional faces of $P$ are polyhedra of dimension $m < n$. *Polygons* (*polytops*) are the polyhedra of dimension $2$ ($3$). A *partition* $\mathcal{P}$ of a polyhedron $P$ is a representation $P = \bigcup_{i=1}^{N} P_i$, where $P_i$ are *subpolyhedra* of $P$ and $\operatorname{int} P_i \cap \operatorname{int} P_j = \emptyset$ if $i \neq j$.

We fix $r \geq 1$ and shall say that something is smooth whenever it is of class $C^r$.

**Definition 2.** (I) *A generalized polyhedron $X$ of dimension $n$ is a closed subset of a smooth manifold $M^n$ and a mapping $f : X \mapsto \mathbf{R}^n$ such that: (1) $f$ extends to a diffeomorphism of an open set $O$, $X \subset O \subset M$, into $\mathbf{R}^n$; (2) $f(X)$ is a polyhedron.*

(II) *A d-dimensional* space $X$ with a generalized polyhedral partition $\mathcal{P}$ *is a subset of a manifold* $M^d$ *and a representation* $X = \bigcup_{i=1}^{n} X_i$ *satisfying the following conditions.*

(i) *The $X_i$ are generalized polyhedra of dimension $d$.*
(ii) $\text{int}(X_i) \cap \text{int}(X_j) = \emptyset$ *if* $i \neq j$.
(iii) *If* $I \subset \{1, \ldots, n\}$ *is such that* $\bigcap_{i \in I} X_i \neq \emptyset$, *the polyhedral structures on* $X_i$, $i \in I$, *agree, so that* $\bigcup_{i \in I} X_i$ *is a generalized polyhedron.*

If $(X, \mathcal{P})$ is as above, we say that $X_i$ are the *atoms* of $\mathcal{P}$ and that $\partial \mathcal{P} = \bigcup_{i=1}^{n} \partial X_i$ is the *boundary* of $\mathcal{P}$.

**Definition 3.** *Let $\mathcal{P}$ and $\mathcal{Q}$ be two partitions of $X$, $\dim X = d$, and let $T : (X, \mathcal{P}) \to (X, \mathcal{Q})$ be such that $T_i = T|_{P_i} : P_i \to Q_i$ and $T_i^{-1} : Q_i \to P_i$ are homeomorphisms, smooth on the interiors. Then we say $T$ is a* generalized polyhedron exchange in $d$ *dimensions. When $d = 2$, we simply speak of a* generalized polygon exchange *(GPE).*

**Examples**. 1. A partition of $X = [0, 1]$ is given by $n$ intervals $P_i = [a_{i-1}, a_i]$, where $0 = a_0 < a_1 < \cdots < a_n = 1$. An interval exchange on the intervals $P_i$, $1 \leq i \leq n$, (see, e.g., [3]) is a generalized polyhedron exchange in one dimension, and the mappings $T_i$ are parallel translations.

2. Let $X$ be a rectangle in $\mathbf{R}^2$, e.g., $X = [0, 1] \times [0, 1]$. Let the atoms of the partition $\mathcal{P} : X = \bigcup_{i=1}^{n} X_i$ be rectangles with the sides parallel to the coordinate axes. Let $t_i$ be $n$ vectors in $\mathbf{R}^2$, such that the rectangles $Q_i = P_i + t_i$, $1 \leq i \leq n$, form a partition, $\mathcal{Q}$, of $X$. This defines a GPE, $T : (X, \mathcal{P}) \to (X, \mathcal{Q})$, where the restrictions, $T_i x = x + t_i$, are translations.

2'. The space $X \subset \mathbf{R}^2$ is an arbitrary polygon; the atoms of a partition $\mathcal{P} : X = \bigcup P_i$, $1 \leq i \leq n$, are subpolygons. We define an *affine polygon exchange* on $\mathcal{P}$ by $n$ affine transformations $T_i : \mathbf{R}^2 \to \mathbf{R}^2$ such that $\mathcal{Q} = \{Q_i = T_i(P_i), i = 1, \ldots, n\}$ is a partition of $X$.

3. An obvious analog of Example 2 in three dimensions features the unit cube $X = [0, 1] \times [0, 1] \times [0, 1] \subset \mathbf{R}^3$ partitioned by $n$ *rectangular parallelepipeds* $P_i$. The mappings $T_i$ are parallel translations $x \mapsto x + t_i$, $t_i \in \mathbf{R}^3$, such that the parallelepipeds $Q_i = P_i + t_i$ form a partition of $X$.

3'. In a three-dimensional version of Example 2', $X$ is a polytop, $\mathcal{P}$ is a partition of $X$ by subpolytops $P_i$, $1 \leq i \leq n$, and $T_i$, $1 \leq i \leq n$, are affine transformations of $\mathbf{R}^3$ such that the polytops $Q_i = T_i(P_i)$, $1 \leq i \leq n$, form a partition, $\mathcal{Q}$, of $X$. This defines an affine polytop exchange $T : (X, \mathcal{P}) \to (X, \mathcal{Q})$.

We say that a partition $\mathcal{P}'$ is a *refinement* of a partition $\mathcal{P}$ (write $\mathcal{P} < \mathcal{P}'$) if every atom of $\mathcal{P}'$ is a subpolyhedron of an atom of $\mathcal{P}$. If $\mathcal{P}$ and $\mathcal{Q}$ are two partitions of $X$, their join $\mathcal{P} \vee \mathcal{Q}$ is the partition formed by the intersections $P \cap Q$ where $P \in \mathcal{P}$ and $Q \in \mathcal{Q}$. For simplicity of exposition, from now on we consider the case $T : (X, \mathcal{R}) \to (X, \mathcal{R})$ and say that $T$ is a *generalized polyhedron exchange on a partition* $\mathcal{R}$. By this we simply mean that $T$ and $T^{-1}$ are continuous on the atoms of $\mathcal{R}$.

3. TOPOLOGICAL ENTROPY

Let $T : X \to X$ be a generalized polyhedron exchange on a partition $\mathcal{R}$. For any finite partition $\mathcal{P}$ we denote by $|\mathcal{P}|$ the number of atoms. We define the $n$th join, $\mathcal{R}_n = \bigvee_{k=0}^{n-1} T^{-k}\mathcal{R}$, so that $T^n : X \to X$ is a generalized polyhedron exchange on a partition $\mathcal{R}_n$. The number $h(T, \mathcal{R}) = \limsup_{n\to\infty} \frac{1}{n} \log |\mathcal{R}_n|$ is called the *entropy* of $T$ *relative to the partition* $\mathcal{R}$, and the *topological entropy* is then given by $h(T) = \sup_{\mathcal{U}} h(T, \mathcal{U})$, where the supremum is over partitions $\mathcal{U}$ such that $\mathcal{R} < \mathcal{U}$. Our definition of the topological entropy for the generalized polyhedron exchanges is similar to the way it is usually defined for homeomorphisms of compacta (see, e.g., [19]). If $\mathcal{R}$ is a *generating partition* (the atoms of the infinite join $\bigvee_{j=0}^{\infty} T^{-j}\mathcal{R}$ consist of single points), then $h(T) = h(T, \mathcal{R})$.

The *growth rate* of a positive sequence $a_n$, $n \geq 1$, is defined by $\vartheta = \limsup_{n\to\infty} \log a_n / n$. We assume that $X$ is endowed with a Finsler metric (i.e., a smooth family of norms, $|| \ ||_x$, on the tangent spaces) and denote by $d(\cdot, \cdot)$ the corresponding distance in $X$. If $\gamma$ is a piecewise $C^1$ curve on $X$, we denote its length by $|\gamma|$. If $\mathcal{P}$ is a partition, we denote by $\partial \mathcal{P}^{(1)}$ the union of edges of $\partial \mathcal{P}$ and set $\ell(\mathcal{P}) = |\partial \mathcal{P}^{(1)}|$. For $x \in X$ let $\mathcal{P}_x \subset \mathcal{P}$ be the set of atoms containing $x$, and set $b(\mathcal{P}) = \max_x |\mathcal{P}_x|$.

**Theorem 4.** *Let $T : X \to X$ be a generalized polyhedron exchange on a partition $\mathcal{R}$, and let $\vartheta$ be the growth rate of the sequence $\ell(\mathcal{R}_n)$. Suppose that the atoms of the partitions $\mathcal{R}_n$, $n > 0$, are connected and that the growth rate of the sequence $b(\mathcal{R}_n)$ is less than or equal to $\vartheta$. Then*

(i) $h(T, \mathcal{R}) \leq \vartheta$;
(ii) *if $\mathcal{R}$ is a generating partition, then $h(T) \leq \vartheta$.*

The bound on the topological entropy given by Theorem 4 is especially useful for the GPEs, for two reasons. First, in two dimensions we can keep track of the growth of the boundary of $\mathcal{R}_n$. Second, we can show that the sequence $b(\mathcal{R}_n)$ does not grow faster than the length, $\ell(\mathcal{R}_n)$.

**Theorem 5.** *Let $T : (X, \mathcal{R}) \to (X, \mathcal{R})$ be a generalized polygon exchange, and let $\vartheta$ be the growth rate of $\ell(\mathcal{R}_n)$. If the atoms of $\mathcal{R}_n$, $n > 0$, are homeomorphic to the disc, then $h(T, \mathcal{R}) \leq \vartheta$. If $\mathcal{R}$ is a generating partition, then $h(T) \leq \vartheta$.*

If the defining partition, $\mathcal{R}$, of a GPE, $T$, is not a generating one, Theorem 5 does not yield a bound on $h(T)$. This is the case, for instance, if $T$ is a diffeomorphism (of a compact surface), since $\mathcal{R}$ is trivial. In this case the *Lyapunov exponent* of $T$ provides an upper bound on $h(T)$ (see, e.g., [13]). Theorem 6 below is an analog of this result in our setting. First, we define the Lyapunov exponent for a GPE.

Let $T : (X, \mathcal{R}) \to (X, \mathcal{R})$ be a GPE. The iterates $T^k$ are smooth on $\text{int}(A)$ for every atom $A$ in the $n$th join $\mathcal{R}_n$, $k \leq n$, and we denote by $DT^k$ the differential of $T^k$. Set $\lambda_k(x) = \log ||DT^k(x)||$, where the norm is with respect to a Finsler metric on $X$. We say $x$ is a *regular point* in $X$, if for any $n > 0$ there exists an atom $A$ in $\mathcal{R}_n$ such that $x \in \text{int}(A)$. In general, the iterates of $T$ are well defined only on the set of regular points, which is dense in $X$.

For regular points in $X$ we define the Lyapunov exponent, $\lambda(x)$, by

$$\lambda(x) = \limsup_{n\to\infty} \frac{\lambda_n(x)}{n} = \limsup_{n\to\infty} \log ||DT^n(x)||^{\frac{1}{n}}.$$

**Theorem 6.** *Let $T : (X, \mathcal{P}) \to (X, \mathcal{P})$ be a generalised polygon exchange, and let $\vartheta \geq 0$ be the growth rate of the sequence $\ell(\mathcal{P}_n)$. Assume that the atoms of $\mathcal{P}_n$ are homeomorphic to the disc. Let $\lambda(x) \leq \vartheta$ for all regular points, and suppose that the following is satisfied*:

*For any $\vartheta' > \vartheta$ there is an integer $N$ such that $\lambda_n(x)/n \leq \vartheta'$ for $n \geq N$ and all regular points $x$.*

*Then the topological entropy of $T$ is bounded above by $\vartheta$.*

*Theorem* 6 *establishes an upper bound on* $h(T)$, *without assuming that $\mathcal{P}$ is a generating partition*; *compare with Theorem* 5.

*Remark* 1. Note that the set, $R_{\mathcal{P}} \subset X$, of regular points for a GPE, depends on the defining partition, $\mathcal{P}$. If $\mathcal{P} < \mathcal{Q}$, then $R_{\mathcal{Q}} \subset R_{\mathcal{P}}$.

*Remark* 2. Let $T : X \to X$ be a GPE on a generating partition, $\mathcal{P}$, and let $\mu$ be a $T$-invariant measure on $X$ supported by $R_{\mathcal{P}}$. We say that $\mu$ is a *regular measure*. Denote by $h_\mu(T)$ the *metric entropy* (see, e.g., [3, 19]). Then $h_\mu(T) \leq h(T, \mathcal{P}) = h(T)$. Thus Theorems 5 and 6 yield an upper bound on the metric entropy of a GPE with respect to a regular measure. We do not know whether the *variational principle* [19] holds for the GPEs.

## 4. APPLICATIONS

1. **Euclidean polygon exchanges.** Let $X \subset \mathbf{R}^2$ be a closed polygon, and let $X = \bigcup_{i=1}^r P_i$ be a polygonal partition. Let $g_i \in \text{Aut}(\mathbf{R}^2)$, $1 \leq i \leq r$, be Euclidean isometries, such that the polygons $Q_i = g_i(P_i)$, $1 \leq i \leq r$, form a partition of $X$. This defines a Euclidean polygon exchange, $T : X \to X$ (see Examples in §2).

**Theorem 7.** *The topological entropy of a Euclidean polygon exchange is zero.*

This is a consequence of Theorem 5. The isometries $g_i$ in the definition of $T$ preserve length. Thus the sequence $|\partial \mathcal{P}_n|$ grows at most linearly, i.e., $\vartheta = 0$.

2. **Polygonal billiards.** Let $\Delta$ be an arbitrary polygon in $\mathbf{R}^2$. We consider the *billiard flow* in $\Delta$ (see, e.g., [3]) and let $T$ be the *billiard ball map*. This means $T$ is the Poincaré map corresponding to the canonical cross-section, $X = X(\Delta) = \partial \Delta \times [0, \pi]$, of the billiard flow. A point $(s, \theta) \in X$ corresponds to a vector with the footpoint $s \in \partial \Delta$, directed inward, whose angle of incidence is $\theta$. Thus $X$ is identified with the rectangle $[0, L] \times [0, \pi]$, where $L = |\partial \Delta|$.

The discontinuity lines of $T$ form a polygonal partition, $\mathcal{P} : X = \bigcup_{i=1}^r P_i$. The boundary of $P_i$ consists of vectors, $v = v(s, \theta)$, in $X$ whose trajectories will run into a vertex of $\Delta$ at the first return.

**Theorem 8.** *The polygonal billiard ball map is a GPE, and $h(T) = 0$.*

To apply Theorem 5, we introduce the Finsler metric $|dm| = |ds|\sin\theta + |d\theta|$ in $X$. The atoms of $\mathcal{P}$ are convex in this metric, and $T$ preserves convexity. Besides, $|\partial \mathcal{P}_n|$ grows at most quadratically; thus $\vartheta = 0$.

*Remark* 3. Equation $h(T) = 0$ implies the subexponential growth rate for the number of vertices in $\mathcal{P}_n$. The latter is equivalent to the main result of [10]. It also implies $h_\mu(T) = 0$ (see [2]) where $d\mu = \sin\theta ds d\theta$ is the standard $T$-invariant measure on $X$.

3. **Billiards in rational polytops.** Let $V \subset \mathbf{R}^3$ be a closed connected polytop. The *billiard flow* in $V$ has a natural cross-section, $X$, and let $T : X \to X$ be the *billiard ball map*. We identify $X$ with $\partial V \times S_+^2$, where $S_+^2 \subset \mathbf{R}^3$ is the upper hemisphere of the unit sphere.

Any collection, $\{L_i : i \in I\}$, of hyperplanes in $\mathbf{R}^3$ defines a subgroup of $O(3)$ generated by the reflections about $L_i$. We denote by $G = G(V)$ the group corresponding to the collection of hyperplanes spanned by the faces, $F_i \subset \partial V$. We say that $V$ is a *rational polytop* if $|G(V)| < \infty$.

Let $V$ be a rational polytop, and let $\Omega \subset S^2$, $\Omega \cong S^2/G$, be a *fundamental polygon* for the action of $G$. The billiard flow in $V$ decomposes into the *directional billiard flows*, $T_\omega^t$, $\omega \in \Omega$. The two-dimensional version of this decomposition is well known [6].

Let $T_\omega : X_\omega \to X_\omega$ be the *directional billiard maps*, where $X = \bigcup_{\omega \in \Omega} X_\omega$ is the decomposition of the standard cross-section. The space $X_\omega$, $\dim X_\omega = 2$, does not depend on $\omega \in \text{int}(\Omega)$, and we denote it by $S = S(V)$. The surface $S$ is a cross-section for the directional billiard flows $T_\omega^t$ realized on the *covering polytop*, $W = W(V)$, glued together from $|G|$ copies of $V$ (see [6] for an analog of this in two dimensions). By construction, $S$ has a partition into $|I||G|/2$ polygons, corresponding to the pairs $(F_i \subset \partial V, g \in G)$.

**Theorem 9.** *Let $V$ be a rational polytop, and let the notation be as above.*

*For any $\omega \in \Omega$ the directional billiard map, $T_\omega : S \to S$, is a GPE, and $h(T_\omega) = 0$.*

Theorem 9 is the three-dimensional version of the well-known fact that the directional billiards in a rational polygon have entropy zero (see, e.g., [3]). It follows from Theorem 5 and an estimate on the growth of the singular set for $T_\omega^n$, as $n \to \infty$.

**Corollary 10.** *Let $T^t$ be the billiard flow in a rational polytop, and let $\mu$ be any invariant measure supported on the set of regular points in the phase space (e.g., the Lebesgue measure, see [3]). Then the metric entropy, $h_\mu(T^t)$, is equal to zero.*

A detailed exposition of the results announced in this note will be published elsewhere [21].

## Acknowledgments

We thank M. Rychlik and N. Chernov for useful discussions.

Mathematics Department, University of Southern California, Los Angeles, California 90089-1113

*E-mail address*, E. Gutkin: `egutkin@math.usc.edu`

*E-mail address*, N. Haydn: `nhaydn@math.usc.edu`